\documentclass[a4paper,10pt]{article}
\usepackage[T1]{fontenc}
\usepackage[utf8]{inputenc}
\usepackage[margin=2.5cm]{geometry}
\usepackage{amsmath,amssymb,amsthm,amscd,mathrsfs}
\usepackage{microtype}
\usepackage[colorlinks=true,linkcolor=blue,citecolor=blue,urlcolor=blue]{hyperref}
\newtheorem{theorem}{Theorem}[section]
\newtheorem{proposition}[theorem]{Proposition}
\newtheorem{lemma}[theorem]{Lemma}

\theoremstyle{definition}
\newtheorem{definition}[theorem]{Definition}
\theoremstyle{remark}
\newtheorem{remark}[theorem]{Remark}
\newtheorem{example}[theorem]{Example}
\newcommand{\A}{\mathbb A}
\newcommand{\Q}{\mathbb Q}
\newcommand{\R}{\mathbb R}
\newcommand{\C}{\mathbb C}
\newcommand{\F}{\mathbb F}
\newcommand{\Z}{\mathbb Z}
\newcommand{\one}{\mathbf 1}
\newcommand{\SB}{\mathscr S}
\DeclareMathOperator{\Tr}{Tr}
\DeclareMathOperator{\Nm}{N}
\DeclareMathOperator{\ord}{ord}
\DeclareMathOperator{\res}{res}

\DeclareMathOperator{\divi}{div}
\DeclareMathOperator{\Spec}{Spec}
\DeclareMathOperator{\Diff}{Diff}

\DeclareMathOperator{\Hom}{Hom}
\title{Harmonic analysis approach to the relative Riemann--Roch theorem on global fields}
\author{Weronika Czerniawska}
\date{}
\begin{document}
\maketitle
\begin{abstract}
In this paper we revisit Tate's harmonic-analytic approach to Riemann--Roch and give a relative formulation in purely adelic terms. For a finite separable extension of global fields, we express the relative Euler characteristic of a divisor as the logarithm of a single adelic integral, using a Haar measure determined by the relative different. Combining this normalization with adelic Poisson summation gives the relative Riemann--Roch formula and the corresponding duality relation. In addition, we define relative ``sizes of cohomology'', extracting the $h^0$ and $h^1$ terms of the Euler characteristic through integrals on the field and its adelic quotient. We make explicit how these numbers depend on the choice of Haar measure. The construction treats the absolute and relative cases together and applies both to arithmetic curves and to nonsingular projective curves over finite fields.
\end{abstract}
\tableofcontents

\section{Introduction}
The work of Tate \cite{tate} and Iwasawa \cite{Iwasawa,Iwasawa_2019} shows how harmonic analysis on adeles brings together the arithmetic of number fields and the geometry of curves over finite fields. For either kind of global field $K$, the adele ring $\A_K$ is a locally compact, Pontryagin self-dual additive group, and the diagonally embedded field $K$ is discrete with compact quotient. These properties allow Poisson summation to take the same form in both settings. In particular, the relation between a divisor and its canonical complement can be expressed through Fourier transformation, giving a harmonic-analytic interpretation of Riemann--Roch and Serre duality; see also \cite{fourier,wei}.

This common description is striking because the invariants arise differently in the two settings. For a curve over a finite field, the Euler characteristic is an alternating sum of dimensions of cohomology groups. For a number field, its arithmetic analogue is expressed through the covolume of a metrized fractional ideal. Adelic harmonic analysis allows both to be studied using the same local functions: characteristic functions of fractional ideals at nonarchimedean places and Gaussians at archimedean places.

The central idea of this paper is to express the Euler characteristic as the logarithm of a single adelic integral and to use the choice of Haar measure to pass to the relative case. Let $K/L$ be a finite separable extension, and let $D_\alpha$ be the divisor represented by an idele $\alpha\in\A_K^\times$. We associate to $\alpha$ a Schwartz--Bruhat function $f_\alpha$ and construct a relative Haar measure $\mu_{K/L}$ for which
\[
\chi_{K/L}(D_\alpha)
=\log\int_{\A_K}f_\alpha(x)\,d\mu_{K/L}(x).
\]
At each nonarchimedean place, the normalization is determined by the local relative different. The integral therefore combines the divisor and the arithmetic of the extension in one expression.

The same test function gives a decomposition of the Euler characteristic into cohomology numbers. Its restriction to $K$ defines a theta sum, while its periodization descends to the quotient $\A_K/K$. After normalizing the periodization to take the value one at the origin, we obtain
\[
h^0_{K/L}(D_\alpha)=\log\sum_{a\in K}f_\alpha(a),
\qquad
h^1_{K/L}(D_\alpha)=-\log\int_{\A_K/K}\bar f_\alpha\,d\mu^{K/L}_{\A_K/K}.
\]
Weil's integration formula then gives
\[
\chi_{K/L}(D_\alpha)=h^0_{K/L}(D_\alpha)-h^1_{K/L}(D_\alpha).
\]
Thus the two cohomology numbers arise from integration on the field and its adelic quotient, which are Pontryagin dual to one another. The family of test functions is preserved by Fourier transformation up to scalar multiplication, and Fourier inversion yields the corresponding duality formula.

The absolute invariants are closely related to the theta invariants of van der Geer and Schoof \cite{schoof} and to Borisov's arithmetic cohomology \cite{borisov}. Our contribution is to give a common adelic integral construction of the Euler characteristic and its two cohomological terms, and to carry this construction through the relative normalization. The numerical Riemann--Roch identities are classical. Their integral realization makes it possible to follow the passage to the relative case at the level of the measures and the functions themselves. In particular, counting measure on $K$ remains fixed, so the change of normalization is reflected entirely in $h^1$. Section~\ref{sec:grr} gives a geometric interpretation of the resulting integral: its logarithm is the degree of the determinant of the pushforward of the line bundle represented by $\alpha$, with the specified metrics in the number-field case.

The dependence on Haar measure is essential to this description. We therefore begin with an arbitrary Haar measure
\[
\mu=c\mu_{0,K},\qquad c>0,
\]
where $\mu_{0,K}$ is the absolute self-dual measure. This allows the absolute and relative formulas to be derived together. For $n=[K:L]$, the relative normalization satisfies
\[
\mu_{K/L}=d_L^{n/2}\mu_{0,K},\qquad
\chi_{K/L}(D)=\chi_K(D)-n\chi_L(0).
\]
In the function-field case, the base term is accounted for by the canonical-divisor identity
\[
K_{X_K}\sim\pi^*K_{X_L}+\operatorname{Diff}(\pi).
\]
The construction includes all places and allows the full constant fields of $K$ and $L$ to differ.

The test functions $f_\alpha$ and their normalized periodizations $\bar f_\alpha$ also belong to the class of cross positive definite functions studied by Banaszczyk \cite{banaszczyk}, whose examples include periodic Gaussians. After the first version of this article was submitted, the author was made aware of Miyada's unpublished work \cite{miyada} on adelic arithmetic cohomology for matrix divisors. His construction considerably overlaps with the ideas developed here. The relation to his work, including the interpretation through Borisov's ghost spaces and the proposed extension to relative vector-bundle invariants, is discussed in Section~\ref{sec:future}.

We state the main formulas below. Throughout the paper, logarithms are natural. Put
\[
d_K=\begin{cases}
|\operatorname{disc}(K/\Q)|,&K\text{ a number field},\\
q_K^{2g_K-2},&K\text{ a function field},
\end{cases}
\]
where $\F_{q_K}$ is the full constant field. Write $\kappa_K$ for the absolute conductor idele associated with the chosen compatible global additive character. The divisor convention and the test functions are specified in Sections~\ref{sec:ideles} and~\ref{sec:general}.

\begin{theorem}[Relative integral formula]\label{thm:relative}
Let $K/L$ be a finite separable extension of global fields of degree $n$. For $\alpha\in\A_K^\times$, the relative measure $\mu_{K/L}=d_L^{n/2}\mu_{0,K}$ gives
\begin{align*}
\chi_{K/L}(D_\alpha)
&=\log\int_{\A_K}f_\alpha\,d\mu_{K/L}\\
&=\log|\alpha|-\frac12\log d_K+\frac n2\log d_L\\
&=\chi_K(D_\alpha)-n\chi_L(0).
\end{align*}
In particular,
\[
\chi_{K/L}(D_\alpha)-\chi_{K/L}(0)=\log|\alpha|.
\]
For number fields, $d_K/d_L^n$ is the numerical relative discriminant. For function fields, it is the norm of the full relative different divisor.
\end{theorem}
This expresses the relative degree and discriminant terms through a single integral. Section~\ref{sec:relative} proves that the local relative normalizations give the stated global measure.

\begin{theorem}[Dependence on Haar measure]\label{thm:scaling}
Fix a compatible global additive character with $K^\perp=K$, and let $\mu=c\mu_{0,K}$, $c>0$. Keep counting measure on $K$, and use the quotient measure determined by Weil's formula. Set
\[
S_\alpha=\sum_{a\in K}f_\alpha(a),\qquad
\bar f_\alpha(\bar u)=S_\alpha^{-1}\sum_{a\in K}f_\alpha(u+a),
\]
and define
\[
h^0_\mu(D_\alpha)=\log S_\alpha,\qquad
h^1_\mu(D_\alpha)=-\log\int_{\A_K/K}\bar f_\alpha\,d\mu_{\A_K/K},
\]
together with
\[
\chi_\mu(D_\alpha)=\log\int_{\A_K}f_\alpha\,d\mu.
\]
Writing a subscript $0$ for the self-dual normalization, one has
\[
h^0_\mu=h^0_0,\qquad h^1_\mu=h^1_0-\log c,\qquad
\chi_\mu=\chi_0+\log c=h^0_\mu-h^1_\mu.
\]
\end{theorem}
These formulas follow from Weil's integration formula and describe how the cohomological decomposition behaves under a change of measure. Their proof is given in Section~\ref{sec:general}.

\begin{theorem}[Duality]\label{thm:duality}
Under the hypotheses of Theorem~\ref{thm:scaling},
\[
h^1_\mu(D_\alpha)=h^0_\mu(D_{\kappa_K\alpha^{-1}})-\log c.
\]
In particular, the self-dual normalization gives absolute duality, and the relative normalization gives
\[
h^1_{K/L}(D_\alpha)=h^0_{K/L}(D_{\kappa_K\alpha^{-1}})
-\frac n2\log d_L.
\]
\end{theorem}
The proof in Section~\ref{sectionserre} uses the Fourier transform of $f_\alpha$ and inversion on $\A_K/K$, retaining the dependence on the quotient measure.

Section~2 establishes the local and global conventions. Section~3 develops the cohomology numbers for arbitrary Haar measure and compares the absolute invariants with their geometric and arithmetic counterparts. Sections~4--6 treat the relative normalization, Poisson summation, and relative duality. Section~\ref{sec:grr} gives the comparison with Grothendieck--Riemann--Roch, and Section~\ref{sec:future} discusses further directions. The local Fourier calculations are collected in Appendix~\ref{apx}.

\subsection*{Acknowledgements}
This work was supported by the SNSF grants 182111 and 204125, as well as by NCCR SwissMAP.

The author thanks the referees for their careful reading, helpful corrections, and valuable suggestions, which have improved the clarity and precision of the manuscript.

\section{Background}
\subsection{Global fields, curves, and discriminants}\label{sec:background}
We treat number fields and function fields of regular projective geometrically integral curves over their full finite fields of constants. In positive characteristic, a global field $K$ is a finite separable extension of $\F_p(t)$ for a suitable separating element $t$.

We define the proper arithmetic curve $X_K$ as follows. For a number field, $X_K$ is the Arakelov compactification of $\Spec O_K$, obtained by adjoining the archimedean places, namely embeddings $K\hookrightarrow\C$ modulo complex conjugation. Here ``proper arithmetic curve'' is used in the Arakelov sense; this compactification is not a proper scheme over $\Spec\Z$. For a function field of characteristic $p$, $X_K$ is the regular projective integral curve over $\F_p$ with function field $K$. The full constant field $\F_{q_K}$ is the algebraic closure of $\F_p$ in $K$. Viewed over $\F_{q_K}$, the curve is geometrically integral. Genus and cohomology dimensions are computed over $\F_{q_K}$.

For number fields, $O_K$ denotes the ring of integers. For a finite separable extension $K/L$ of function fields, choose the same separating element $t\in L$ and let $O_L,O_K$ be the integral closures of $\F_p[t]$ in $L,K$. These affine rings omit the points above the pole of $t$. Define the codifferent and different by
\[
\mathfrak D_{K/L}^{-1}=\{x\in K:\Tr_{K/L}(xO_K)\subseteq O_L\},
\qquad \mathfrak D_{K/L}=(\mathfrak D_{K/L}^{-1})^{-1}.
\]
The different is an ideal of $O_K$. The relative discriminant $\mathfrak d_{K/L}$ is the ideal of $O_L$ whose localization at a nonzero prime $\mathfrak p$ is generated by
\[
\det\bigl(\Tr_{K/L}(b_i b_j)\bigr)_{1\leq i,j\leq n},\qquad n=[K:L],
\]
where the $b_i$ form an $O_{L,\mathfrak p}$-basis of $O_K\otimes_{O_L}O_{L,\mathfrak p}$. A global integral basis is not required. These local ideals are independent of the bases and satisfy $\mathfrak d_{K/L}=\Nm_{K/L}(\mathfrak D_{K/L})$. Their supports are the ramified primes in the respective rings. The numerical norm
\[
d_{K/L}=\#(O_L/\mathfrak d_{K/L})
\]
is a positive integer, not an ideal. Standard facts on differents and discriminants may be found in \cite[Chapter III, Section 2]{Neu} and \cite{serre}.

Separately, set
\begin{equation}\label{eq:absolute-d}
d_K=\begin{cases}|\operatorname{disc}(K/\Q)|,&K\text{ a number field},\\
q_K^{2g_K-2},&K\text{ a function field}.
\end{cases}
\end{equation}
The function-field invariant is genus-based: it is not in general the norm of the different of $K/\F_p(t)$. For example, $d_{\F_p(t)}=p^{-2}$, whereas the different of the identity extension is trivial.

For the finite separable morphism $\pi:X_K\to X_L$ of function-field curves, put
\[
\Diff(\pi)=\sum_v\delta_{v/w}[v],\qquad
\mathfrak D_{K_v/L_w}=\pi_v^{\delta_{v/w}}O_v.
\]
The sum includes \emph{all} closed points. We use the distinct notation
\begin{equation}\label{eq:full-different}
\Delta_{K/L}=\prod_v\#k(v)^{\delta_{v/w}}
\quad\text{in the function-field case},\qquad
\Delta_{K/L}=d_{K/L}\quad\text{in the number-field case}.
\end{equation}
Thus in positive characteristic
\[
\Delta_{K/L}=d_{K/L}\prod_{v\text{ above a pole of }t}\#k(v)^{\delta_{v/w}}.
\]
The affine norm equals the full norm if $K/L$ is unramified at all those omitted points. We make no such assumption and retain their contributions. This avoids imposing an extra choice of $t$ for the measure construction.

\subsection{Places, local moduli, and duality}
\begin{definition}
Two nontrivial absolute values on $K$ are equivalent if they induce the same topology. A place is an equivalence class of nontrivial absolute values.
\end{definition}
Equivalent absolute values are positive powers of one another. The identity on $K$ extends uniquely to a topological-field isomorphism of their completions. The resulting topological field $K_v$ depends only on the place; its specified absolute value still depends on the representative. If $|a|_2=|a|_1^b$ on $K$, $b>0$, the same relation holds on the identified completions, and the canonical isomorphism need not be an isometry.

Every $K_v$ is locally compact. For a Haar measure $\mu_v$, define the normalized local modulus by
\begin{equation}\label{haarnorm}
|a|_v=\frac{\mu_v(aV)}{\mu_v(V)}\quad(a\ne0),\qquad |0|_v=0,
\end{equation}
where $V$ is measurable and $0<\mu_v(V)<\infty$. The ratio is independent of $V$ and of Haar normalization. At real places it is the usual absolute value; at complex places it is its square. The latter need not satisfy the triangle inequality and is called a valuation in the broader convention of \cite{cassels}. We use ``modulus'' to distinguish it from the absolute values defining places.

For a nonarchimedean place let $O_v$ be the valuation ring, $\pi_v$ a uniformizer, $k(v)$ the residue field, and $\nu_v:K_v^\times\to\Z$ the discrete valuation. Then
\[
|a|_v=\#k(v)^{-\nu_v(a)}\quad(a\in K_v^\times).
\]
For $v\mid w$ and $a\in K_v^\times$, one has $|a|_v=|\Nm_{K_v/L_w}(a)|_w$. The product formula is
\[
\prod_v|a|_v=1\qquad(a\in K^\times).
\]

For a locally compact abelian group $G$, define its Pontryagin dual
\[
\widehat G=\Hom_{\mathrm{cont}}(G,S^1),\qquad S^1=\{z\in\C:|z|=1\},
\]
with pointwise multiplication and the compact-open topology. The evaluation map identifies $G$ with $\widehat{\widehat G}$. A group is Pontryagin self-dual if it admits a topological-group isomorphism with its dual. A choice of nontrivial additive character $\psi_v$ gives such an isomorphism
\[
K_v\longrightarrow\widehat{K_v},\qquad a\longmapsto(x\longmapsto\psi_v(ax)).
\]
This depends on the character and is distinct from self-duality of a measure.

For this fixed pairing, write $\mu_{0,v}$ for the self-dual measure. If $\mu_v=b\mu_{0,v}$, its dual, transported to $K_v$, is $b^{-1}\mu_{0,v}$. The two are equal precisely when $b=1$. In particular, the trace pairing on an unramified extension of $\Q_p$ has self-dual measure with $\mu_{0,v}(O_v)=1$; self-duality is not restricted to $\Q_p$.

\subsection{Adeles and compatible global characters}\label{sec:characters}
The adele ring is the restricted product
\[
\A_K=\prod_v'K_v
=\{(a_v):a_v\in O_v\text{ for almost all nonarchimedean }v\},
\]
with its restricted product topology. It is a locally compact ring, and $K\subset\A_K$ is discrete with compact quotient.

Let $f_v\in\SB(K_v)$, with $f_v=\one_{O_v}$ at almost all nonarchimedean places. Their restricted tensor product $f=\bigotimes_v'f_v\in\SB(\A_K)$ is the function
$f((a_v))=\prod_v f_v(a_v)$. At each adele only finitely many factors differ from $1$. For measures with $\mu_v(O_v)=1$ almost everywhere,
\begin{equation}\label{eq:restricted-integral}
\int_{\A_K}f\,d\mu=\prod_v\int_{K_v}f_v\,d\mu_v.
\end{equation}
Indeed, choose a finite set $S$ containing the archimedean places and every exception for the functions and measures. The function vanishes outside $\prod_{v\in S}K_v\times\prod_{v\notin S}O_v$; the compact tail has measure $1$. Finite product integration on the remaining factors proves the formula. It is not a claim about arbitrary infinite products of functions.

We now fix globally compatible characters. For number fields use on $\Q_p$ and $\R$, respectively,
\[
\psi_{0,p}(x)=\exp(2\pi i\{x\}_p),\qquad
\psi_{0,\infty}(x)=\exp(-2\pi ix),
\quad \{x\}_p=\sum_{j=-k}^{-1}a_jp^j.
\]
Here $x=\sum_{j\ge-k}a_jp^j$ is its $p$-adic expansion. For $v\mid w$ define $\psi_v=\psi_{0,w}\circ\Tr_{K_v/\Q_w}$. For $a\in\Q$, $a-\sum_p\{a\}_p\in\Z$, so the product of the base characters on $a$ is $1$. Compatibility of local and global traces gives the same assertion on $K$.

For a function field choose a nonzero rational differential $\omega$ and set
\begin{equation}\label{eq:residue-character}
\psi_v(x)=\exp\left(\frac{2\pi i}{p}
\Tr_{k(v)/\F_p}\res_v(x\omega)\right).
\end{equation}
An element of $\F_p$ in the exponent is represented by any integer lift. The residue theorem, followed by trace to $\F_p$, shows that $\psi=\prod_v\psi_v$ is trivial on $K$. For comparisons in a finite separable extension we may choose a nonzero $\omega_L$ and take $\omega_K=\pi^*\omega_L$. This pullback is nonzero because the extension is separable.

The global self-duality theorem identifies
\begin{equation}\label{eq:annihilator}
\A_K\xrightarrow{\sim}\widehat{\A_K},\quad
 a\longmapsto(x\longmapsto\psi(ax)),\qquad K^\perp=K.
\end{equation}
We use here the full annihilator theorem, not merely triviality of $\psi$ on $K$. For number fields see \cite[Sections 3.2 and 4.1]{tate} and \cite[Section 4, Functional equation]{kudla}; the global-field statement with residue pairings is part of adelic duality in \cite{wei}. One can also reduce the annihilator statement to $\Q$ or $\F_p(t)$ by the nondegenerate trace pairing for a finite separable extension; multiplication of the differential by an element of $K^\times$ does not change $K^\perp$. Consequently
\[
\widehat{\A_K/K}=K,\qquad \widehat K=\A_K/K.
\]
An arbitrary set of nontrivial local characters would not suffice for these identifications.

Our Fourier convention, throughout, is
\begin{equation}\label{eq:fourier}
\widehat f_\mu(y)=\int_{\A_K}f(x)\psi(xy)\,d\mu(x),\qquad
f(x)=\int_{\A_K}\widehat f_\mu(y)\psi(-xy)\,d\widehat\mu(y).
\end{equation}
Both formulas apply to Schwartz--Bruhat functions. With $\mu=c\mu_{0,K}$, $\widehat\mu=c^{-1}\mu_{0,K}$ and
\begin{equation}\label{eq:fourier-square}
\mathcal F_\mu=c\mathcal F_0,\qquad
\mathcal F_\mu^2f(x)=c^2f(-x).
\end{equation}
The second identity uses the same measure twice; inversion using the dual measure has no extra factor.

\subsection{Ideles and canonical divisors}\label{sec:ideles}
The idele group $\A_K^\times$ consists of $(\alpha_v)$ with $\alpha_v\in K_v^\times$ everywhere and $\alpha_v\in O_v^\times$ almost everywhere. Its modulus is $|\alpha|=\prod_v|\alpha_v|_v$. Our divisor convention is
\[
D_\alpha=\sum_{v\nmid\infty}-\nu_v(\alpha_v)[v]
+\sum_{v\mid\infty}\log|\alpha_v|_v[v]
\]
for number fields, with Arakelov degree assigning $\log\#k(v)$ to a finite point and degree $1$ to each archimedean symbol. Thus $\widehat{\deg}D_\alpha=\log|\alpha|$. For function fields the divisor is $\sum_v-\nu_v(\alpha_v)[v]$ and
$\log|\alpha|=(\deg_{\F_{q_K}}D_\alpha)\log q_K$.
Multiplication by $a\in K^\times$ changes $D_\alpha$ by a principal divisor, with this sign convention.

Define the absolute conductor idele by
\begin{equation}\label{eq:kappa}
O_v^\perp=\kappa_{K,v}O_v\quad(v\text{ nonarchimedean}),\qquad
\kappa_{K,v}=1\quad(v\text{ archimedean}).
\end{equation}
For number fields $O_v^\perp=\mathfrak D_{K_v/\Q_p}^{-1}$. For function fields write $\omega=h_v\,dt_v$ in a local parameter; then $O_v^\perp=h_v^{-1}O_v$. These statements are proved in Appendix~\ref{apx}. Hence $D_{\kappa_K}$ is the absolute canonical divisor in our conventions; in positive characteristic it is $\divi(\omega)$. In both cases,
\begin{equation}\label{eq:kappa-norm}
|\kappa_K|=d_K,\qquad \mu_{0,v}(O_v)=|\kappa_{K,v}|_v^{-1/2}.
\end{equation}
At archimedean places the self-dual measures are Lebesgue measure on $\R$ and twice Lebesgue measure on $\C$.

The relative idele $\kappa_{K/L}$ is defined separately by local generators of $\mathfrak D_{K_v/L_w}^{-1}$, with archimedean components $1$. Its divisor is the relative different divisor, and $|\kappa_{K/L}|=\Delta_{K/L}$. For function fields the local differential formula gives
\begin{equation}\label{eq:canonical}
\divi(\pi^*\omega_L)=\pi^*\divi(\omega_L)+\Diff(\pi).
\end{equation}
Equivalently, $K_{X_K}\sim\pi^*K_{X_L}+\Diff(\pi)$. The number-field counterpart is transitivity of the different, with the stated trivial archimedean metrics. A relative different must not be identified with the absolute canonical divisor without the base term.

\subsection{Quotient measures and periodization}\label{submu}
A function $F:\A_K\to\C$ is $K$-periodic if $F(u+a)=F(u)$ for every $u\in\A_K$ and $a\in K$. It descends uniquely to $\A_K/K$. The periodization of $f\in\SB(\A_K)$ is the particular $K$-periodic function
\[
Pf(\bar u)=\sum_{a\in K}f(u+a).
\]
For the test functions used below the sum converges locally uniformly, and hence gives a continuous quotient function. In number fields the finite components restrict the sum, locally in $u$, to finitely many translates of fractional-ideal lattices in $K\otimes_\Q\R$; Gaussian decay gives uniform convergence on compact subsets. In function fields the functions have compact support, and discreteness of $K$ makes the relevant sums locally finite. The same argument applies to the Fourier transforms computed in Appendix~\ref{apx}. These hypotheses justify periodization, termwise integration, and the Poisson formulas below; compare \cite[Section 4.2]{tate}.

Keep counting measure $\mu_K$ on $K$. For $\mu$ on $\A_K$, Weil's formula determines $\mu_{\A_K/K}$ by
\begin{equation}\label{weil}
\int_{\A_K}f\,d\mu
=\int_{\A_K/K}\sum_{a\in K}f(u+a)\,d\mu_{\A_K/K}(\bar u).
\end{equation}
For a measurable fundamental domain $B$ and a continuous $K$-periodic $F$ this gives
\begin{equation}\label{intd}
\int_{\A_K/K}F\,d\mu_{\A_K/K}=\int_B F\,d\mu.
\end{equation}
Indeed, restriction to a fundamental domain realizes the quotient integral; the translates of $B$ partition $\A_K$ up to null sets, and summing their integrals proves \eqref{weil}. See also \cite[Chapter 5, Section 11.1]{Hav}.

For the compatible self-dual measure, $\mu_{0,\A_K/K}(\A_K/K)=1$. To recall why normalization matters, let this volume initially be $C$. Poisson summation for a lattice and its annihilator gives $\sum_K f=C^{-1}\sum_K\widehat f_0$. Applying it again to $\widehat f_0$, using $K^\perp=K$ and $\mathcal F_0^2 f(x)=f(-x)$, gives $C^2=1$ for a test function with nonzero sum. Since $C>0$, $C=1$. For an arbitrary rescaling the factor in \eqref{eq:fourier-square} compensates for the changed volume, so the same reasoning does not force that volume to be $1$.

Thus, if $\mu=c\mu_{0,K}$,
\begin{equation}\label{eq:quotient-volume}
\mu_{\A_K/K}=c\mu_{0,\A_K/K},\qquad
\mu_{\A_K/K}(\A_K/K)=c.
\end{equation}
The measure dual to $\mu_{\A_K/K}$ on $\widehat{\A_K/K}=K$ has mass $c^{-1}$ at each point. In contrast, the measure on $\widehat K=\A_K/K$ dual to counting measure is probability measure. It equals the quotient measure only if $c=1$. This distinction is exactly the reciprocal basic-domain-volume factor discussed after \cite[Theorem 4.2.1, pp.~333--334]{tate}.

\section{Adelic cohomology numbers for arbitrary Haar measure}\label{sec:general}
Fix the compatible character of Section~\ref{sec:characters}. Let $\mu_0=\mu_{0,K}$ and $\mu=c\mu_0$, $c>0$. Counting measure on $K$ and the quotient measure in \eqref{weil} remain fixed by these choices. For an idele $\alpha$, set $f_\alpha(x)=f_1(\alpha^{-1}x)$, where
\begin{equation}\label{eq:test-functions}
f_{\alpha_v}(x)=\begin{cases}
\one_{O_v}(\alpha_v^{-1}x),&v\text{ nonarchimedean},\\
\exp\bigl(-e_v\pi|\alpha_v^{-1}x|_v^{2/e_v}\bigr),&v\text{ archimedean},
\end{cases}
\qquad f_\alpha=\bigotimes_v'f_{\alpha_v},
\end{equation}
and $e_v=[K_v:\R]$. These are Schwartz--Bruhat functions and are integrable; the restricted product argument is \eqref{eq:restricted-integral}, in accord with \cite[Lemma 3.3.2]{tate}. As mentioned before, they form a Fourier-stable family up to scalar multiplication.

Their theta sums and normalized periodizations are
\begin{equation}\label{eq:theta}
S_\alpha=\sum_{a\in K}f_\alpha(a),\qquad
\bar f_\alpha(\bar u)=S_\alpha^{-1}Pf_\alpha(\bar u).
\end{equation}
Convergence was explained in Section~\ref{submu}; $S_\alpha\ge1$ since $f_\alpha(0)=1$. In particular, $\bar f_\alpha$ is well defined and $\bar f_\alpha(\bar0)=1$.

\begin{definition}\label{def:cohomology}
The logarithmic cohomology numbers and Euler characteristic for $\mu$ are
\begin{align}
h^0_\mu(D_\alpha)&=\log\int_K f_\alpha\,d\mu_K=\log S_\alpha,\label{eq:h0}\\
h^1_\mu(D_\alpha)&=-\log\left(\int_{\A_K/K}\bar f_\alpha\,d\mu_{\A_K/K}\right),\label{eq:h1}\\
\chi_\mu(D_\alpha)&=\log\left(\int_{\A_K}f_\alpha\,d\mu\right).\label{eq:chi}
\end{align}
\end{definition}
Units at finite places and phases at complex places do not change these functions. For $a\in K^\times$, changing $\alpha$ to $a\alpha$ reindexes the theta sum, while the product formula preserves the adelic integral. Weil's formula then preserves the quotient integral. Thus the numbers depend only on the divisor class (with its archimedean metric in the number-field case).

The real and complex Gaussian integrals with the chosen measures equal $1$. Equations~\eqref{eq:restricted-integral} and~\eqref{eq:kappa-norm} therefore give
\begin{equation}\label{eq:base-integral}
\int_{\A_K}f_1\,d\mu_0=d_K^{-1/2}.
\end{equation}
In function fields this is $q_K^{1-g_K}$, the product of the conductor contributions of the rational differential, including infinity.
By change of variables and Weil's formula,
\begin{align}
\int_{\A_K}f_\alpha\,d\mu&=c|\alpha|d_K^{-1/2},\label{eq:integral-general}\\
\int_{\A_K/K}\bar f_\alpha\,d\mu_{\A_K/K}
&=\frac{c|\alpha|d_K^{-1/2}}{S_\alpha}.\label{eq:period-integral}
\end{align}
Taking logarithms proves
\begin{equation}\label{eq:RR-general}
\chi_\mu(D_\alpha)=h^0_\mu(D_\alpha)-h^1_\mu(D_\alpha)
=\log|\alpha|-\tfrac12\log d_K+\log c.
\end{equation}
It also proves Theorem~\ref{thm:scaling}: the theta sum is unchanged, the quotient integral is multiplied by $c$, and its negative logarithm changes by $-\log c$.

\subsection{Fourier duality}\label{sectionserre}
The local calculations in Appendix~\ref{apx} yield
\begin{equation}\label{eq:global-transform}
\widehat f_{\alpha,\mu}
=c|\alpha|d_K^{-1/2}f_{\kappa_K\alpha^{-1}}.
\end{equation}
The absolute conductor $\kappa_K$ is independent of the normalization of $\mu$. Define the quotient Fourier coefficient by integration against $\psi(ua)$ with the quotient measure. Weil's formula and $a\in K=K^\perp$ give
\begin{equation}\label{eq:quotient-transform}
\widehat{\bar f_\alpha}_{\mu_{\A_K/K}}(a)
=S_\alpha^{-1}\widehat f_{\alpha,\mu}(a)\qquad(a\in K).
\end{equation}
This is the normalized version of the periodization identity in \cite[Lemma 4.2.3]{tate}. It involves a Fourier coefficient, not the value of $\bar f_\alpha$ at $a$.

Fourier inversion on $\A_K/K$ uses the dual discrete measure of mass $c^{-1}$ per point. Evaluating at zero, where $\bar f_\alpha(0)=1$, gives
\[
S_\alpha=\frac1c\sum_{a\in K}\widehat f_{\alpha,\mu}(a)
=|\alpha|d_K^{-1/2}S_{\kappa_K\alpha^{-1}}.
\]
The Fourier series is absolutely convergent by \eqref{eq:global-transform} and convergence of the corresponding theta sum, so this inversion is justified. Substitution in \eqref{eq:period-integral} yields
\begin{equation}\label{eq:dual-integral}
\int_{\A_K/K}\bar f_\alpha\,d\mu_{\A_K/K}
=\frac{c}{S_{\kappa_K\alpha^{-1}}}.
\end{equation}
Its negative logarithm is the formula in Theorem~\ref{thm:duality}. Equivalently, inversion reads
\[
1=c^{-1}S_{\kappa_K\alpha^{-1}}
\int_{\A_K/K}\bar f_\alpha\,d\mu_{\A_K/K}.
\]
The factor $c^{-1}$ belongs inside the logarithm; it does not multiply $h^0$ or $h^1$.

\subsection{Absolute invariants and the geometric interpretation}\label{sec:comparison}
Write $h_K^i=h_0^i$ and $\chi_K=\chi_0$. In a function field,
\[
S_\alpha=\#H^0(X_K,\mathcal O(D_\alpha))
=q_K^{\dim_{\F_{q_K}}H^0(X_K,\mathcal O(D_\alpha))}.
\]
Indeed the finite-place membership conditions are precisely
$\ord_v(a)+\operatorname{coeff}_v(D_\alpha)\ge0$. Absolute duality and the usual geometric Serre duality give
\begin{equation}\label{eq:geometric}
h_K^i(D)=\dim_{\F_{q_K}}H^i(X_K,\mathcal O(D))\log q_K,
\qquad i=0,1,
\end{equation}
and
\[
\chi_K(D)=(\deg_{\F_{q_K}}D+1-g_K)\log q_K.
\]
Thus the natural-logarithm invariants are logarithmic cardinalities, not literally the vector-space dimensions. Equivalently, the harmonic-analytic identity gives the numerical Riemann--Roch relation between the spaces of sections of $D$ and $K_{X_K}-D$.

For a number field put
\[
I_\alpha=\{a\in K:a\in\alpha_vO_v\text{ for every finite }v\},\qquad
\|a\|_\alpha^2=\sum_{v\mid\infty}e_v|\alpha_v^{-1}a|_v^{2/e_v}.
\]
Then
\[
S_\alpha=\sum_{a\in I_\alpha}\exp(-\pi\|a\|_\alpha^2).
\]
This is the theta sum defining the effectivity invariant in \cite{schoof}. The lattice $I_\alpha$ with Gaussian $a\mapsto\exp(-\pi\|a\|_\alpha^2)$ also encodes the data used for Borisov's first-kind ghost space \cite[Section 5]{borisov}. Our absolute $h^0$ is the logarithm of that same sum, while $h^1$ is its dual-lattice counterpart and $\chi_K=\log|\alpha|-\frac12\log d_K$ is the usual arithmetic Euler characteristic. This identifies the numerical ingredients, without claiming that our relative construction provides additional cohomology spaces.

One can recover the absolute normalization from any $\mu=c\mu_0$ by dividing by its quotient volume:
\begin{align*}
\chi_K(D_\alpha)&=\log\left(
\frac{\int_{\A_K}f_\alpha\,d\mu}{\mu_{\A_K/K}(\A_K/K)}\right),\\
h_K^1(D_\alpha)&=-\log\left(
\frac{\int_{\A_K/K}\bar f_\alpha\,d\mu_{\A_K/K}}
{\mu_{\A_K/K}(\A_K/K)}\right).
\end{align*}
We retain the raw measure-dependent integrals in Definition~\ref{def:cohomology}; these normalized expressions are a distinct convention.

\section{Relative Euler characteristic}\label{sec:relative}
Let $K/L$ be a finite separable extension of degree $n$. At every nonarchimedean place $v\mid w$, define the local relative measure by
\begin{equation}\label{eq:relative-local}
\mu_{K_v/L_w}(O_v)=\#k(v)^{-\delta_{v/w}/2}
=\Nm_{L_w}\bigl(\mathfrak d_{K_v/L_w}\bigr)^{-1/2}.
\end{equation}
The local discriminant is an ideal of $O_{L_w}$, and its local norm is the index of that ideal. This is the inverse square root of its numerical norm. At archimedean places keep Lebesgue measure on $\R$ and twice Lebesgue measure on $\C$. For function fields every place, including those above the poles of $t$, is included in \eqref{eq:relative-local}. Almost all volumes equal $1$, so the restricted product measure $\mu_{K/L}$ exists and
\begin{equation}\label{eq:relative-base}
\int_{\A_K}f_1\,d\mu_{K/L}=\Delta_{K/L}^{-1/2}.
\end{equation}

We verify the global scale. For number fields the discriminant tower formula is
\begin{equation}\label{disc}
d_K=d_L^n d_{K/L},\qquad
\log d_{K/L}=\log d_K-n\log d_L.
\end{equation}
For function fields the same numerical relation with the full different is
\begin{equation}\label{eq:disc-full}
\Delta_{K/L}=d_K/d_L^n.
\end{equation}
To see explicitly the constant-field factors, put $r=[\F_{q_K}:\F_{q_L}]$. For any divisor $E$ on $X_L$, the local identity $\sum_{v\mid w}e(v/w)f(v/w)=n$ gives
\[
\deg_{\F_{q_K}}(\pi^*E)\log q_K
=n\deg_{\F_{q_L}}(E)\log q_L.
\]
Here $\log q_K=r\log q_L$; the degree formula is therefore $\deg_{\F_{q_K}}\pi^*E=(n/r)\deg_{\F_{q_L}}E$. Taking the weighted degrees in \eqref{eq:canonical} gives
\[
(2g_K-2)\log q_K
=n(2g_L-2)\log q_L+\sum_v\delta_{v/w}\log\#k(v),
\]
which proves \eqref{eq:disc-full}. This retains the canonical contribution of the base in the discriminant quotient formula; compare \cite{det}. If one uses only the affine discriminant, the additional product at infinity in Section~\ref{sec:background} must be included, or unramifiedness at all omitted places must be imposed explicitly.

By uniqueness of Haar measure, comparison of \eqref{eq:relative-base} with \eqref{eq:base-integral} now proves
\begin{equation}\label{eq:relative-scale}
\mu_{K/L}=c\mu_{0,K},\qquad c=d_L^{n/2}.
\end{equation}
For function fields this is $c=q_L^{n(g_L-1)}$. Thus the local-product construction agrees with the scale in Theorem~\ref{thm:relative}, in both settings. The equality is one of global measures; at individual places the absolute differential also contributes the pulled-back base conductor.

Changing variables by $\alpha$ proves Theorem~\ref{thm:relative}, with the additional expression
\begin{equation}\label{eq:relative-chi}
\chi_{K/L}(D_\alpha)=\log|\alpha|-\tfrac12\log\Delta_{K/L}.
\end{equation}
In particular its constant term is
\[
-\tfrac12\log\Delta_{K/L}
=-\tfrac12\log d_K+\tfrac n2\log d_L.
\]
The well-defined comparison as a ratio of integrals is
\begin{equation}\label{eq:ratio}
\chi_{K/L}(D_\alpha)=\log\left(
\frac{\displaystyle\int_{\A_K}f_\alpha\,d\mu_{0,K}}
{\left(\displaystyle\int_{\A_L}f_1\,d\mu_{0,L}\right)^n}\right).
\end{equation}
Indeed the denominator equals $d_L^{-n/2}=c^{-1}$, which proves the formula without introducing any additional function on $\A_K$ pulled from $\A_L$.

\begin{remark}[The base field correction]\label{rem:base}
For number fields $d_\Q=1$ and $\chi_\Q(0)=0$, so $\chi_{K/\Q}=\chi_K$. For $L=\F_p(t)$, however, $d_L=p^{-2}$ and $\chi_L(0)=\log p$. Consequently
\[
\chi_K(D)=\chi_{K/\F_p(t)}(D)+[K:\F_p(t)]\log p.
\]
The absolute geometric theory is therefore not obtained by simply identifying it with the relative normalization over $\F_p(t)$.
\end{remark}

\section{Riemann--Roch and Poisson summation}\label{sec:poisson}
For a fixed measure, change of variables alone gives
\[
\chi_\mu(D_\alpha)-\chi_\mu(0)=\log|\alpha|.
\]
Together with the theta-sum duality of Section~\ref{sectionserre}, this is the integral formulation of classical numerical Riemann--Roch. For number fields it has the usual arithmetic degree and discriminant terms, as in \cite[Chapter V, Section 1]{lang} and \cite[Chapter III]{Neu}; for function fields \eqref{eq:geometric} gives the customary degree-and-genus formula. The relative formula is its normalization in \eqref{eq:relative-scale}.

We spell out the scaling in Poisson summation to distinguish the relative Euler correction from the invariant theta identity. Put $g_\alpha(x)=f_1(\alpha x)$. For the relative measure,
\[
\widehat g_{\alpha,\mu_{K/L}}(y)
=|\alpha|^{-1}\widehat f_{1,\mu_{K/L}}(\alpha^{-1}y)
=\frac{\Delta_{K/L}^{-1/2}}{|\alpha|}
 f_1(\kappa_K^{-1}\alpha^{-1}y).
\]
Since the quotient volume is $c=d_L^{n/2}$, Poisson summation gives
\begin{equation}\label{eq:relative-poisson}
\sum_{a\in K}f_1(\alpha a)
=\frac{\Delta_{K/L}^{-1/2}}{c|\alpha|}
\sum_{a\in K}f_1(\kappa_K^{-1}\alpha^{-1}a).
\end{equation}
Here $\Delta_{K/L}^{-1/2}/c=d_K^{-1/2}$. Thus \eqref{eq:relative-poisson} is exactly the absolute Poisson identity. Rescaling the Haar measure changes the transform and the reciprocal quotient-volume factor in opposite ways. The relative correction to the Euler characteristic enters through the integral defining $h^1$, rather than a new theta-sum identity.

Section~\ref{sec:grr} describes the adelic restriction of scalars for $K/L$ and compares the relative integral with the determinant degree supplied by Grothendieck--Riemann--Roch.

\section{Relative cohomology numbers}\label{sec:relative-cohomology}
The relative cohomology numbers are the specialization of Definition~\ref{def:cohomology} to $\mu_{K/L}$. Denote the quotient measure by $\mu^{K/L}_{\A_K/K}$. It satisfies
\[
\int_{\A_K}f\,d\mu_{K/L}
=\int_{\A_K/K}\sum_{a\in K}f(u+a)\,d\mu^{K/L}_{\A_K/K}(\bar u).
\]
Counting measure $\mu_K$ is still used on $K$; the dual discrete measure, which assigns mass $c^{-1}$ at each point, is different in general.

\begin{definition}
For $\alpha\in\A_K^\times$, define
\begin{align*}
h^0_{K/L}(D_\alpha)&=\log\int_K f_\alpha\,d\mu_K=\log S_\alpha,\\
h^1_{K/L}(D_\alpha)&=-\log\left(
\int_{\A_K/K}\bar f_\alpha\,d\mu^{K/L}_{\A_K/K}\right).
\end{align*}
\end{definition}
The signs agree with \eqref{eq:h0}--\eqref{eq:h1}; there is a single minus sign outside the second logarithm and no additional inverse. Directly,
\[
\int_{\A_K/K}\bar f_\alpha\,d\mu^{K/L}_{\A_K/K}
=S_\alpha^{-1}\int_{\A_K}f_\alpha\,d\mu_{K/L},
\]
so $h^0_{K/L}-h^1_{K/L}=\chi_{K/L}$. Equations~\eqref{eq:quotient-volume} and~\eqref{eq:relative-scale} give
\begin{align*}
\mu^{K/L}_{\A_K/K}(\A_K/K)&=c=d_L^{n/2},\\
h^0_{K/L}(D_\alpha)&=h_K^0(D_\alpha),\\
h^1_{K/L}(D_\alpha)&=h_K^1(D_\alpha)-\tfrac n2\log d_L.
\end{align*}
The dependence is additive. Only the exponentiated quantities satisfy a multiplicative formula, namely $\exp(h^1_{K/L})=c^{-1}\exp(h_K^1)$. These relative logarithmic quantities need not be nonnegative, and should not be interpreted as dimensions of vector spaces.

\subsection{Relative duality}\label{relser}
Theorem~\ref{thm:duality} specializes to
\begin{equation}\label{eq:relative-serre}
h^1_{K/L}(D_\alpha)
=h^0_{K/L}(D_{\kappa_K\alpha^{-1}})-\tfrac n2\log d_L.
\end{equation}
The Fourier proof uses the absolute $\kappa_K$, even though the measure is relative. Inversion at zero, with mass $c^{-1}$ on the dual discrete group, reads
\[
1=c^{-1}S_{\kappa_K\alpha^{-1}}
\int_{\A_K/K}\bar f_\alpha\,d\mu^{K/L}_{\A_K/K}.
\]
Taking logarithms proves \eqref{eq:relative-serre}. This also explains why using counting measure on the dual discrete group without changing the quotient measure would give the wrong formula.

\begin{example}[The rational function field]\label{ex:rational}
Let $K=L=\F_p(t)$, $\alpha=1$, and choose $\omega=dt$. At infinity put $s=t^{-1}$; then $dt=-s^{-2}ds$, so $\kappa_{K,\infty}O_\infty=s^2O_\infty$ and $\mu_{0,\infty}(O_\infty)=p$. All finite $O_v$ have absolute volume $1$. Hence $\int f_1\,d\mu_0=p$ and $S_1=p$. The absolute quotient has volume $1$, the normalized-periodization integral is $1$, and
\[
h_K^0(0)=\log p,\qquad h_K^1(0)=0,\qquad \chi_K(0)=\log p.
\]
For the identity relative extension the local differents are trivial, all local relative volumes are $1$, and $c=d_K^{1/2}=p^{-1}$. Thus the quotient volume and the normalized-periodization integral are both $p^{-1}$, giving
\[
h^0_{K/K}(0)=\log p,\qquad h^1_{K/K}(0)=\log p,\qquad
\chi_{K/K}(0)=0.
\]
More generally, for the relative measures over $\F_p(t)$ the quotient volume is $p^{-[K:\F_p(t)]}$, whereas the absolute self-dual quotient volume is $1$.
\end{example}

\begin{example}[A constant-field extension]
For $L=\F_q(t)$ and $K=\F_{q^r}(t)$, one has $n=r$, $d_K=q^{-2r}$, $d_L=q^{-2}$, and $\Delta_{K/L}=1$. Therefore $c=q^{-r}$ and $\chi_{K/L}(0)=0$, while $\chi_K(0)=r\log q$. This checks both the constant-field degree factor and the base correction in an unramified extension.
\end{example}

\section{Comparison with Grothendieck--Riemann--Roch}\label{sec:grr}
Let $K/L$ be a finite separable extension of degree $n$. We explain how the adelic data used in this paper represent the objects entering Riemann--Roch, and how the relative integral computes the degree of the determinant of a pushforward. Throughout this section, $\deg_F$ denotes $\deg_{\F_{q_F}}(\cdot)\log q_F$ for function-field curves and the unnormalized arithmetic degree for number fields.

\paragraph{From ideles to pushforward.}
For a function field $F$, the standard identification
\[
\operatorname{Pic}(X_F)
\simeq F^\times\backslash\A_F^\times/\prod_vO_v^\times
\]
sends the class of $\alpha$ to $\mathcal L_\alpha=\mathcal O_{X_F}(D_\alpha)$. Recall our convention $D_\alpha=\sum_v-\nu_v(\alpha_v)[v]$. The local lattices $\alpha_vO_v$ therefore describe $\mathcal L_\alpha$: for $a\ne0$, membership is equivalent to $\nu_v(a)-\nu_v(\alpha_v)\geq0$, while $a=0$ is automatically included.

For a number field, the corresponding presentation is
\[
\widehat{\operatorname{Pic}}(O_F)
\simeq F^\times\backslash\A_F^\times/
\left(\prod_{v\nmid\infty}O_v^\times\times\prod_{v\mid\infty}U_v\right),
\qquad U_v=\{u\in F_v^\times:|u|_v=1\},
\]
where $\widehat{\operatorname{Pic}}(O_F)$ is the group of isometry classes of metrized line bundles. The local lattices $\alpha_vO_v$ at the nonarchimedean places determine the underlying fractional ideal $I_\alpha$, while the archimedean components specify the fibre norms
\[
\|z\|_\tau=|\alpha_v|_v^{-1/e_v}|z|
\]
at embeddings $\tau:F\hookrightarrow\C$ belonging to $v$. Thus the indicators and Gaussians defining $f_\alpha$ encode the line bundle and its metric.

Let $\pi:X_K\to X_L$ be the finite morphism of curves in the geometric case, and $\pi:\Spec O_K\to\Spec O_L$ the finite morphism in the arithmetic case. Its pushforward is restriction of scalars. At a nonarchimedean place $w$ of $L$, the local lattice of $E_\alpha=\pi_*\mathcal L_\alpha$ is
\[
\prod_{v\mid w}\alpha_vO_v
\subseteq K\otimes_LL_w\simeq\prod_{v\mid w}K_v.
\]
The product $\prod_{v\mid w}f_{\alpha_v}$ is its indicator. At archimedean places, the same grouping of Gaussian factors represents the orthogonal-sum metric on the pushforward. This only regroups the factors of $f_\alpha$.

An $L$-basis of $K$ induces $T:\A_L^n\xrightarrow{\sim}\A_K$ with $T(L^n)=K$. Both absolute self-dual measures give quotient volume one, so, with $\nu_L=\mu_{L/L}=d_L^{1/2}\mu_{0,L}$,
\begin{equation}\label{eq:grr-measure}
T_*(\mu_{0,L}^{\otimes n})=\mu_{0,K},
\qquad T_*(\nu_L^{\otimes n})=\mu_{K/L}.
\end{equation}
The second equality follows from \eqref{eq:relative-scale}. These identities are independent of the basis by the product formula.

\paragraph{The geometric comparison.}
In the function-field case, view both curves over $\F_{q_L}$. The line bundle represented by $\alpha$ determines a class in the Grothendieck group $K_0(X_K)$ of vector bundles; see \cite[Chapter II, Section 8]{weibelK}. Since $\pi$ is finite, higher direct images vanish and $\pi_![\mathcal L_\alpha]=[E_\alpha]$.

Write $\operatorname{CH}^*(X)_\Q=\operatorname{CH}^*(X)\otimes\Q$. The relative form of the Grothendieck--Riemann--Roch diagram is
\begin{equation}\label{eq:grr-diagram}
\begin{CD}
K_0(X_K)\otimes\Q
 @>{\operatorname{ch}(-)\,\operatorname{Td}_\pi}>> \operatorname{CH}^*(X_K)_\Q\\
@V{\pi_!}VV @VV{\pi_*}V\\
K_0(X_L)\otimes\Q @>{\operatorname{ch}}>> \operatorname{CH}^*(X_L)_\Q.
\end{CD}
\end{equation}
Here $\operatorname{Td}_\pi=\operatorname{Td}(T_{X_K})/\pi^*\operatorname{Td}(T_{X_L})$, with tangent bundles taken over $\F_{q_L}$. For curves,
\[
\operatorname{ch}(\mathcal L_\alpha)=1+[D_\alpha],
\qquad \operatorname{Td}_\pi=1-\frac12[\operatorname{Diff}(\pi)],
\]
where the second identity follows from \eqref{eq:canonical}. The degree-one part of Grothendieck--Riemann--Roch gives
\begin{equation}\label{eq:grr-degree-class}
c_1(E_\alpha)=\pi_*\left([D_\alpha]-\frac12[\operatorname{Diff}(\pi)]\right).
\end{equation}
Taking logarithmic degrees yields
\begin{align*}
\deg_L\det E_\alpha
&=\deg_KD_\alpha-\frac12\deg_KK_{X_K}
 +\frac n2\deg_LK_{X_L}\\
&=\log|\alpha|-\frac12\log\Delta_{K/L}\\
&=\chi_K(D_\alpha)-n\chi_L(0).
\end{align*}
This convention retains the formula when the full constant fields differ: if their extension degree is $r$, restriction of scalars multiplies cohomology dimensions by $r$, while $\log q_K=r\log q_L$.

\paragraph{The arithmetic comparison.}
For number fields, put $\overline E_\alpha=\pi_*\overline{\mathcal L}_\alpha$, with the orthogonal-sum metrics described above. The determinant and norm formulas in \cite[Chapter III, (7.2)--(7.4)]{Neu} give
\[
\det\overline E_\alpha
\simeq N_{K/L}(\overline{\mathcal L}_\alpha)
\otimes\det(\pi_*\overline{\mathcal O}_K).
\]
Here $N_{K/L}$ is the norm of metrized line bundles: the ideal norm at finite places and the product metric on the tensor product of fibres over embeddings extending each embedding of $L$. Moreover,
\[
\bigl(\det\pi_*\overline{\mathcal O}_K\bigr)^{\otimes2}
\simeq\overline{\mathfrak d}_{K/L},
\]
where the discriminant ideal is equipped with the trivial metric. Hence
\[
\deg_L\det\overline E_\alpha
=\log|\alpha|-\frac12\log d_{K/L}
=\chi_K(D_\alpha)-n\chi_L(0).
\]

In both settings, comparison with \eqref{eq:relative-chi} gives, with the specified metrics understood in the arithmetic case,
\begin{equation}\label{eq:grr-integral}
\deg_L\det E_\alpha
=\log\int_{\A_K}f_\alpha\,d\mu_{K/L}
=\chi_{K/L}(D_\alpha).
\end{equation}
Thus the integral is the exponential of the determinant degree. Its Haar normalization supplies the numerical different correction appearing in Riemann--Roch. Weil's formula expresses this logarithm as $h^0_{K/L}-h^1_{K/L}$; the comparison identifies their alternating sum, without identifying the separate relative numbers with ordinary sheaf-cohomology dimensions.

\section{Future perspectives}\label{sec:future}
The first natural generalization of this work is the treatment of vector bundles. The author has been working on an extension of the theory developed in this paper, inspired by \cite[Section 1]{morishita}. After the first version of this work was submitted, the author was made aware of the unpublished work of Miyada \cite{miyada}, which overlaps considerably with the ideas in this article. Miyada \cite{miyada} already treats the vector-bundle case through matrix divisors over global fields, constructing adelic cohomology as ghost spaces and establishing Borisov--Serre duality and Riemann--Roch. Miyada's work and our planned follow-up on vector bundles appear to share an inspiration in both Morishita's treatment of matrix divisors and adelic test functions in \cite[Section 1]{morishita} and Borisov's ghost spaces \cite{borisov}. The distinction relevant here is that Miyada does not develop the relative case for a finite extension of global fields, nor does he explicitly formulate the Euler characteristic as the logarithm of a single adelic integral, although the integral computation underlying that identity occurs in his proof of Theorem~3.

In rank one, our absolute $h^1$ should agree with the dimension of Miyada's ghost space of the second kind on $\A_K/K$. Its defining probability measure is obtained by normalizing the periodization of $f_\alpha$ to have total mass one. The relative $h^1$ is expected to be the dimension of the same ghost space with respect to the relative quotient Haar measure, in Borisov's measure-dependent sense \cite[Definition 2.4]{borisov}. The probability measure and the underlying convolution structure remain unchanged, while the change of reference Haar measure changes the density at zero and hence the dimension. Since counting measure in the definition of $h^0$ is unchanged, the relative correction to the Euler characteristic is carried entirely by $h^1$. Therefore, it is expected that, by extending Miyada's ideas, we can interpret our relative cohomology numbers as dimensions of concrete adelic ghost spaces.

It would also be interesting to apply the numerical cohomology theory of Price \cite{num} in the adelic setting and to the relative, rather than only the absolute, case. This would require comparing its axioms and constructions with the measure-dependent logarithmic numbers defined above.

It is also natural to investigate the relation of this framework to the theory of adelic curves developed by Chen and Moriwaki \cite{chen}.

A further question is whether these methods can contribute to an adelic formulation of the Arakelov adjunction formula or of Faltings--Riemann--Roch for arithmetic surfaces. Beilinson adeles and their cohomological interpretation \cite{hub}, together with completed adelic constructions on arithmetic surfaces \cite{wc}, suggest possible settings. These are prospective applications; the one-dimensional arguments given here do not establish them.

\appendix
\section{Fourier analysis on local fields}\label{apx}
We supply the local statements used in the global transform. All Fourier transforms have kernel $\psi_v(xy)$; inversion uses $\psi_v(-xy)$ and the dual measure, as in \eqref{eq:fourier}.

\subsection{Characters and conductors}
For $\Q_p$, the fractional-part expression is
\[
\{x\}_p=\sum_{j=-k}^{-1}a_jp^j,\qquad
x=\sum_{j\ge-k}a_jp^j,\quad 0\le a_j<p.
\]
The character is $\exp(2\pi i\{x\}_p)$; the real character is $\exp(-2\pi ix)$. These opposite signs ensure global compatibility. For a finite extension $E/\Q_p$, use $\psi_E=\psi_{0,p}\circ\Tr_{E/\Q_p}$. The definition of the codifferent gives
\[
O_E^\perp=\{y:\psi_E(xy)=1\text{ for all }x\in O_E\}
=\mathfrak D_{E/\Q_p}^{-1}.
\]
Indeed the kernel of the base character is $\Z_p$, so the condition is precisely $\Tr(yO_E)\subseteq\Z_p$.

In positive characteristic use \eqref{eq:residue-character}. Writing $\omega=h_v\,dt_v$, the annihilator is $h_v^{-1}O_v$. If $yh_v\in O_v$ then $xyh_v\,dt_v$ has zero residue for every $x\in O_v$. Conversely, if $yh_v$ has a negative Laurent term, multiplication by an appropriate nonnegative power of $t_v$ and a residue-field coefficient produces a residue with nonzero trace to $\F_p$. Nondegeneracy of the finite-field trace pairing proves the assertion. At infinity on $\F_p(t)$, the calculation $dt=-s^{-2}ds$ is therefore essential; a coefficient-of-$t^{-1}$ rule chosen independently at every place would miss this contribution.

\subsection{Orthogonality and local transforms}
\begin{lemma}\label{L1}
Let $E/\Q_p$ be finite, $\pi$ a uniformizer, $q=\#(O_E/\pi O_E)$, and $\mathfrak D_{E/\Q_p}=\pi^\delta O_E$. For any Haar measure $\mu$ and any integer $m$,
\[
\int_{\pi^mO_E}\psi_E(x)\,d\mu(x)
=\begin{cases}q^{-m}\mu(O_E),&m\ge-\delta,\\0,&m<-\delta.\end{cases}
\]
\end{lemma}
\begin{proof}
The character is trivial on $\pi^mO_E$ exactly when that ideal is contained in $\mathfrak D^{-1}=\pi^{-\delta}O_E$. In that case the integral is its volume. Otherwise choose $z$ in the subgroup with $\psi_E(z)\ne1$. Translation by $z$ multiplies the integral by $\psi_E(z)$ and also leaves it unchanged, so it vanishes.
\end{proof}

\begin{lemma}\label{L2}
With the same notation,
\[
\widehat{\one_{\pi^mO_E}}_\mu(y)
=q^{-m}\mu(O_E)\one_{\pi^{-m-\delta}O_E}(y).
\]
More generally, for any nonarchimedean local field with the chosen pairing and $O_v^\perp=\kappa_{K,v}O_v$,
\begin{equation}\label{eq:local-transform}
\widehat{\one_{\alpha_vO_v}}_{\mu_v}(y)
=|\alpha_v|_v\mu_v(O_v)
\one_{\alpha_v^{-1}\kappa_{K,v}O_v}(y).
\end{equation}
\end{lemma}
\begin{proof}
The character $x\mapsto\psi_v(xy)$ is trivial on $\alpha_vO_v$ exactly when $y\in\alpha_v^{-1}O_v^\perp$. Orthogonality as in Lemma~\ref{L1} gives the volume in this case and zero otherwise. The volume is $|\alpha_v|_v\mu_v(O_v)$.
\end{proof}
Applying this formula twice to $\one_{O_v}$ with the self-dual measure shows
$\mu_{0,v}(O_v)^2|\kappa_{K,v}|_v=1$, which proves \eqref{eq:kappa-norm}.

\begin{example}[A ramified check]
If $p$ is odd and $E=\Q_p(\sqrt p)$, then $\pi=\sqrt p$ is a uniformizer and the different is generated by $2\pi$, hence equals $\pi O_E$. The trace character is trivial on $\pi^{-1}O_E$. The integral in Lemma~\ref{L1} at $m=-1$ is therefore its positive volume. This illustrates why a cutoff at $m=0$ would be incorrect at ramified places.
\end{example}

\subsection{Archimedean factors and the global formula}
At a real place the usual Gaussian transform with Lebesgue measure gives
$\widehat f_{\alpha_v}(y)=|\alpha_v|_v f_{\alpha_v^{-1}}(y)$.
At a complex place write $z=x+iy$ and $w=u+iv$. The character is
$\exp(-2\pi i\Tr_{\C/\R}(zw))=\exp(-4\pi i(xu-yv))$, and the measure is $2\,dx\,dy$. The two real Gaussian integrals give
\[
\mathcal F\bigl(e^{-2\pi|z|^2}\bigr)(w)=e^{-2\pi|w|^2}.
\]
After scaling, the same formula $\widehat f_{\alpha_v}=|\alpha_v|_v f_{\alpha_v^{-1}}$ follows, where the modulus is the squared usual absolute value. Since $\kappa_{K,v}=1$, this is the archimedean instance of the family transformation. If the local measure is rescaled by $b$, the right-hand side is multiplied by $b$.

\begin{proposition}\label{a4}
For any Haar measure $\mu=c\mu_{0,K}$ on $\A_K$ and any idele $\alpha$,
\begin{align*}
\widehat f_{\alpha,\mu}&=\widehat f_{\alpha,\mu}(0)
 f_{\kappa_K\alpha^{-1}},\\
\widehat f_{\alpha,\mu}(0)&=|\alpha|\int_{\A_K}f_1\,d\mu
=c|\alpha|d_K^{-1/2}.
\end{align*}
For the relative measure the coefficient is $|\alpha|\Delta_{K/L}^{-1/2}$; it equals $|\alpha|d_{K/L}^{-1/2}$ for number fields and for function-field affine presentations unramified at every omitted point.
\end{proposition}
\begin{proof}
Multiply the local formulas of Lemma~\ref{L2} and the archimedean Gaussian formulas using restricted product integration. Their supports give $f_{\kappa_K\alpha^{-1}}$, and their scalar factors give the integral of $f_\alpha$. Changing variables gives $|\alpha|\int f_1\,d\mu$, whose value is \eqref{eq:integral-general}. For the relative measure use \eqref{eq:relative-base}.
\end{proof}
Thus changing the measure changes the scalar, while the absolute conductor remains fixed. In particular, the factor $|\alpha|$ must be retained. With $\alpha=1$ one also has the unambiguous evaluation identity
\[
\widehat f_{1,\mu_{K/L}}(\alpha^{-1}u)
=\Delta_{K/L}^{-1/2}f_1(\kappa_K^{-1}\alpha^{-1}u),
\]
which is the expression used in Section~\ref{sec:poisson}.

\paragraph{Data availability statement}
All mathematical arguments and calculations supporting this study are contained in the paper. An earlier version is available at \href{https://arxiv.org/abs/2208.10424}{arXiv:2208.10424}.

% Bibliography is embedded below so the native single-file editor can compile it.
% The accompanying ad.bib preserves the supplied BibTeX database and adds the
% references requested by the referees and added in subsequent revisions.

\bigskip
\noindent W.~Czerniawska, \textsc{Westlake University, China}\\
\textit{E-mail address}: \texttt{weronika.czerniawska@westlake.edu.cn}
\end{document}